\documentclass{article}
\usepackage{amsmath,amssymb,enumerate}

\newcommand{\nosymbol}{}
\newcommand{\tmem}[1]{{\em #1\/}}
\newcommand{\tmop}[1]{\ensuremath{\operatorname{#1}}}
\newcommand{\um}{-}
\newenvironment{enumerateroman}{\begin{enumerate}[i.] }{\end{enumerate}}
\newenvironment{itemizeminus}{\begin{itemize} }{\end{itemize}}
\newtheorem{lemma}{Lemma}
\newtheorem{theorem}{Theorem}

\begin{document}

\title{An Easton-like Theorem for Zermelo-Fraenkel Set Theory Without
Choice\\
(Preliminary Report)}\author{\and Anne Fernengel and Peter Koepke\\
Mathematical Institute, University of Bonn, Germany}\maketitle

In Zermelo-Fraenkel set theory without the Axiom of Choice the
$\theta$-function
\[ \theta \left( \kappa \right) = \sup \left\{ \nu \in \text{Ord} \left| \right.  \text{there
   is a surjection } f : \mathcal{P} \left( \kappa \right) \rightarrow \nu
   \right\} \]
provides a ``surjective'' substitute for the exponential function
$2^{\kappa}$. It generalizes the cardinal $\theta = \theta \left( \aleph_0
\right)$ which is prominent in descriptive set theory. Obviously
\begin{itemizeminus}
  \item $\theta \left( \kappa \right)$ is a cardinal $> \kappa^+ $;
  
  \item $\kappa < \lambda$ implies $\theta \left( \kappa \right) \leqslant
  \theta \left( \lambda \right)$.
\end{itemizeminus}
We show that these are the only restrictions that $\tmop{ZF}$ imposes on
$\theta(\kappa)$.

\begin{theorem}
  Let $M$ be a ground model of $\tmop{ZFC} + \tmop{GCH} +${\tmem{Global
  Choice}}. In $M$, let $F$ be a function defined on the class of infinite
  cardinals such that
  \begin{enumerateroman}
    \item $F \left( \kappa \right)$ is a cardinal $> \kappa^+$;
    
    \item $\kappa < \lambda$ implies $F \left( \kappa \right) \leqslant F
    \left( \lambda \right) $.
  \end{enumerateroman}
  Then there is an extension $N$ of $M$ which satisfies $\tmop{ZF}$, preserves
  cardinals and cofinalities, and such that $\theta \left( \kappa \right) = F
  \left( \kappa \right)$ holds for all cardinals in $N$.
\end{theorem}

This is a version of Easton's Theorem \cite{Easton} for {\tmem{all}} infinite cardinals,
irrespective of cofinalities.

The model $N$ is defined as an inner model of a class generic extension $M
\left[ G \right]$ of $M$. In this note we only indicate main elements of the
construction. Some techniques are inspired by \cite{Gitik-Koepke}. Detailed proofs and further investigations of the model and its
variants, in particular concerning cardinal arithmetic properties and the
amount of choice possible, will be presented in the PhD thesis of the first
author.

Let
\[ A = \bigcup_{\kappa \in \tmop{Card}} \left\{ \kappa \right\} \nosymbol
   \times F \left( \kappa \right) \]
We say that $t = \left( t, <_t \right)$ is an $F$-{\tmem{tree}} if
\begin{itemizeminus}
  \item $t \subseteq A$ and $t$ is the field of the binary relation $<_t $;
  
  \item $<_t$ is strict and transitive;
  
  \item $\left( \kappa, \nu \right) <_t \left( \lambda, \mu \right)
  \rightarrow \kappa < \lambda$;
  
  \item $\left( \lambda, \mu \right) \in t \wedge \kappa \in \tmop{Card} \cap
  \lambda \rightarrow \exists ! \nu \left( \kappa, \nu \right) <_t \left(
  \lambda, \mu \right)$; hence the predecessors of $\left( \lambda, \mu
  \right) \in t$ with $\lambda = \aleph_{\alpha}$ are linearly ordered by
  $<_t$ in ordertype $\alpha$, i.e., $t$ is a tree;
  
  \item there are finitely many maximal elements $\left( \lambda_0, \mu_0
  \right), \ldots, \left( \lambda_{n - 1}, \mu_{n - 1} \right) \in t$ such
  that
  \[ t = \left\{ \left( \kappa, \nu \right)  \left| \right. \exists i < n
     \left( \kappa, \nu \right) \leqslant_t \left( \lambda_i, \mu_i \right)
     \right\} . \]
\end{itemizeminus}
Let the class forcing $P = \left( P, \leqslant \right)$ consist of all
conditions
\[ p : \left( t, <_t \right) \rightarrow V \]
such that $\left( t, <_t \right)$ is an $F$-tree, and for all $\left( \lambda,
\mu \right) \in t$:
\begin{itemizeminus}
  \item if $\lambda = \aleph_0$ then $p \left( \lambda, \mu \right) \in
  \tmop{Fn} \left( \aleph_0, 2, \aleph_0 \right)$
  
  \item if $\lambda = \kappa^+$ is the successor of $\kappa \in
  \tmop{Card}$ then $p \left( \lambda, \mu \right) \in \tmop{Fn} \left( \left[
  \kappa, \lambda \right), 2, \lambda \right)$
  
  \item otherwise $p \left( \lambda, \mu \right) = \emptyset$
\end{itemizeminus}
Here $\tmop{Fn} \left( D, 2, \lambda \right)$ is the Cohen forcing\\ $\left\{ h
\left| \right. h : \tmop{dom} \left( h \right) \rightarrow 2, \tmop{dom}
\left( h \right) \subseteq D, \tmop{card} \left( h \right) < \lambda
\right\}$, partially ordered by $\supseteq$.

For $p, q \in P$ with $p : \left( t, <_t \right) \rightarrow V$ and $q :
\left( s, <_s \right) \rightarrow V$ set $p \leqslant q$ iff
\begin{itemizeminus}
  \item $t \supseteq s$ and $<_t \supseteq <_s$
  
  \item $\forall \left( \lambda, \mu \right) \in s \text{ } p \left( \lambda,
  \mu \right) \supseteq q \left( \lambda, \mu \right)$
\end{itemizeminus}

Let $G$ be $M$-generic for $P$. For $\left( \lambda, \mu \right) \in A$ let
\[ c_{\lambda, \mu} = \left\{ \xi \left| \right. \exists p \in G
   \left( p : \left( t, <_t \right) \rightarrow V \wedge \exists \left( \kappa, \nu \right) \leqslant_t \left(
   \lambda, \mu \right) p \left( \kappa, \nu \right) \left( \xi \right) = 1
   \right) \right\} \]
be the ``$\mu$-th Cohen subset'' of $\lambda$.

For limit $\delta < F \left( \lambda \right)$ define surjections
\[ S_{\lambda, \delta} : \left\{ c_{\lambda, \mu}  \left|  \right. \mu <
   \delta \right\} \rightarrow \left\{ 0 \right\} \cup \left( \tmop{Lim} \cap
   \delta \right) \]
by
\[ S_{\lambda, \delta} \left( c_{\lambda, \mu} \right) = \mu^{\ast} \]
where $\mu^{\ast}$ is the largest element of $\left\{ 0 \right\} \cup \left(
\tmop{Lim} \cap \delta \right)$ which is $\leqslant \mu$.

Around each $S_{\lambda, \delta}$ define a ``cloud'' $\tilde{S}_{\lambda,
\delta}$ of similar functions:

$\tilde{S}_{\lambda, \delta} = \{ S_{\lambda, \delta} \circ \sigma
\left| \right. \sigma \text{ is a permutation of } \delta \text{, which is the
identity except}$

$\text{for finitely many intervals of the form } \left[ \gamma,
\gamma + \omega \right) \}$.

The predicate $\mathcal{S}$ collects the equivalence classes:
\[ \mathcal{S}= \bigcup \left\{ \tilde{S}_{\lambda, \delta}  \left| \right.
   \lambda \in \tmop{Card} \text{ and } \delta \text{ is a limit ordinal } < F
   \left( \lambda \right) \right\} . \]

The model $N$ is generated over $M$ by the Cohen sets, the surjections
$S_{\lambda, \delta}$ and the predicate $\mathcal{S}$.
\[ N = \left( \tmop{HOD}_{\mathcal{S}} \left( M \cup \left\{ c_{\lambda, \mu} 
   \left| \right.  \left( \lambda, \mu \right) \in A \right\} \cup \left\{
   S_{\lambda, \delta}  \left| \right. \left(\lambda, \delta\right)
   \in A \wedge \delta \in \tmop{Lim} \right\} \right) \right)^{M
   \left[ G \right]} \]
The analysis of $N$ is based on symmetry properties of $P$.

An $F$-{\tmem{permutation}} is a bijection $\pi : A \leftrightarrow A$ such
that
\[ \forall \left( \lambda, \mu \right) \in A \exists \mu' < F \left( \lambda
   \right) \pi \left( \lambda, \mu \right) = \left( \lambda, \mu' \right), \]
i.e., $\pi$ preserves the levels of $A$.

An $F$-permutation $\pi$ canonically extends to a bijection of the class of
all $F$-trees and to an automorphism of $P$.

An $F$-permutation $\pi$ is called {\tmem{small}} if for all assignments $\pi
\left( \lambda, \mu \right) = \left( \lambda, \mu' \right)$ there is $\gamma
\in \left\{ 0 \right\} \cup \tmop{Lim}$ such that $\mu, \mu' \in \left[
\gamma, \gamma + \omega \right)$.

\begin{lemma}
  If $p, q \leqslant r$ and $p \upharpoonright \tmop{dom} \left( r \right)$,
  $q \upharpoonright \tmop{dom} \left( r \right)$ are compatible then there
  exists a small $F$-permutation $\pi$ such that
  \begin{itemizeminus}
    \item $\pi \upharpoonright \tmop{dom} \left( r \right) = \tmop{id}$ and so
    $\pi \left( r \right) = r$
    
    \item $p$ and $\pi \left( q \right)$ are compatible in $P$.
  \end{itemizeminus}
\end{lemma}

A small $F$-permutation will only permute the Cohen sets $\dot{c}_{\lambda,
\mu}^H$ by a finite difference in the second index $\mu$. Since the functions
$S_{\lambda, \delta}$ are invariant with respect to such small pertubations,
one obtains:

\begin{lemma}
  Let $\varphi$ be an $\in$-formula, $\left( \kappa_0, \mu_0 \right), \ldots,
  \left( \kappa_{m - 1}, \mu_{m - 1} \right) \in A$, and\\ $\left( \lambda_0,
  \delta_0 \right), \ldots, \left( \lambda_{n \um 1}, \delta_{n - 1} \right)
  \in A$ with limit ordinals $\delta_i $. Let $\vec{x} \in M$. Let $p : \left(
  t, <_t \right) \rightarrow V$ be a condition such that $\left( \kappa_0,
  \mu_0 \right), \ldots, \left( \kappa_{m - 1}, \mu_{m - 1} \right) \in t$.
  Let $\left( \bar{t}, <_t \right) \subseteq \left( t, <_t \right)$ be the
  subtree $\bar{t} = \left\{ \left( \kappa, \nu \right)  \left| \right.
  \exists i < m \left( \kappa, \lambda \right) \leqslant_t \left( \kappa_i,
  \mu_i \right) \right\}$. Then, using canonical names,
  \[ p \Vdash \varphi \left( \dot{c}_{\kappa_0, \mu_0}, \ldots,
     \dot{c}_{\kappa_{m - 1}, \mu_{m - 1}}, \overrightarrow{\check{x}},
     \dot{S}_{\lambda_0, \delta_0}, \ldots, \dot{S}_{\lambda_{n - 1},
     \delta_{n - 1}}, \dot{\mathcal{S}} \right) \]
  iff
  \[ p \upharpoonright \bar{t} \Vdash \varphi \left( \dot{c}_{\kappa_0,
     \mu_0}, \ldots, \dot{c}_{\kappa_{m - 1}, \mu_{m - 1}},
     \overrightarrow{\check{x}}, \dot{S}_{\lambda_0, \delta_0}, \ldots,
     \dot{S}_{\lambda_{n - 1}, \delta_{n - 1}}, \dot{\mathcal{S}} \right) . \]
\end{lemma}

This implies the following approximation property:

\begin{lemma}
  Let $X \in N$, $X \subseteq \tmop{Ord}$. Then there are indices\\ $\left(
  \kappa_0, \mu_0 \right), \ldots, \left( \kappa_{m - 1}, \mu_{m - 1} \right)
  \in A$ such that $X \in M \left[ c_{\kappa_0, \mu_0}, \ldots, c_{\kappa_{m -
  1}, \mu_{m - 1}} \right]$.
\end{lemma}

By the lemma, $N$ is a cardinal preserving extension of $M$. For all cardinals
$\lambda \in M$ and limit ordinals $\delta < F \left( \lambda \right)$,
$S_{\lambda, \delta} \in N$ yields a surjection from $\mathcal{P} \left(
\lambda \right) \cap N$ onto $\left( \left\{ 0 \right\} \cup \tmop{Lim}
\right) \cap \delta$, thus $\theta \left( \lambda \right) \geqslant F \left(
\lambda \right)$. Further permutation arguments show that $N \models \theta
\left( \lambda \right) = F \left( \lambda \right)$, as required.

\end{document}